\documentclass [12pt,leqno] {article}

\usepackage{amsfonts}

\def\dist{\mathop{\rm dist}\nolimits}
\def\supp{\mathop{\rm supp}\nolimits}

\newcommand{\qed}{~\hfill~$\fbox{}$ }
\newcommand{\proces}{( X_t,P^x)} 

\newcommand{\R}{ \mathbb{R}^{d}}

\newcommand{\V}{\mathbb{V}}
\newcommand{\dowod}{{\em Proof}.\/ }
\newcommand{\indyk}[1]{{\bf 1}_{#1}}
\newcommand{\sfera}{ \mathbb{S}}

\newcommand{\Fourier}{ {\cal F}}
\newcommand{\Dynkin}[2]{ {\cal U} #1 (#2) }
\newcommand{\Dynkinr}[2]{ {\cal U}_r #1 (#2) }
\newcommand{\nusmall}[1]{\hat{\nu}_{#1}}
\newcommand{\gener}{{\cal A}}
 \def\dist{\mathop{\rm
    dist}\nolimits} 

\newtheorem{lemat}{\indent\sc Lemma}
\newtheorem{prop}{\indent\sc Proposition}
\newtheorem{twierdzenie}{\indent\sc Theorem}
\newtheorem{wniosek}[lemat]{\indent\sc Corollary}
\newtheorem{definicja}{\indent\sc Definition}

\newcounter{conum} 

\begin{document}

\title{Estimates of potential kernel and
Harnack's  inequality
for anisotropic fractional Laplacian}
\author{Krzysztof Bogdan \and Pawe{\l} Sztonyk}
\date{June 14, 2005}
\maketitle

\begin{center}
  Abstract
\end{center}
\begin{scriptsize}
We characterize those homogeneous translation invariant symmetric non-local operators 
with positive maximum principle whose harmonic functions 
satisfy Harnack's  inequality.
We also estimate the corresponding semigroup and 
the potential kernel.
\end{scriptsize}

\footnotetext{2000 {\it MS Classification}:
Primary 47D03, 31C05; Secondary 60J35, 60G51.
Research partially supported by KBN and RTN
(HPRN-CT-2001-00273-HARP) }

\section{Main results and background}\label{s:m}
\setcounter{equation}{0}
Let $\alpha\in (0,2)$ and $d\in \{1,2,\ldots\}$.
We consider an arbitrary L\'evy measure on $\R\setminus\{0\}$ which is
symmetric, homogeneous: $\nu(rB)=r^{-\alpha}\nu(B)$, and
nondegenerate (for definitions see Section~\ref{Prel}).
$\nu$ 
yields a convolution semigroup of probability 
measures $\{P_t\,,\;t> 0\}$ on $\R$. 
Each $P_t$ has a smooth density $p_t$. 
We consider the corresponding {\it potential measure} 
$\V=\int_0^\infty P_t \,dt$
and the {\it potential kernel} 
$$  V(x)=\int_0^\infty p_t(x)dt\,,\quad x\in \R\,.
$$
$V(x)=|x|^{\alpha-d}V(x/|x|)$, but it may be infinite in some
directions (\cite[pp. 148-149]{BS}). It is of interest
to study continuity of $V$ on the unit sphere $\sfera$ in $\R$
under specific assumptions on $\nu$ (see (\ref{nu_gamma})).
\begin{twierdzenie}\label{t:V}
If $d>\alpha$ and $\nu$ is a $\gamma$-measure on $\sfera$ with
$\gamma>d-2\alpha$ then $V$ is continuous on $\sfera$.
\end{twierdzenie}
The following partial converse
shows that the threshold $d-2\alpha$ is exact.
\begin{twierdzenie}
  \label{t:nu}
If $\V$ is a $\kappa$-measure on $\sfera$ then $\nu$ is a
$(\kappa-2\alpha)$-measure on $\sfera$.
\end{twierdzenie}
In particular, if $V$ is bounded on $\sfera$ then
$\nu$ is a $(d-2\alpha)$--measure on $\sfera$.

We define an operator $\gener$ on smooth 
functions $\varphi$ with compact support in $\R$, $\varphi \in C^\infty_c(\R)$, by
\begin{eqnarray}
  \gener\varphi(x) &=&
  \int\limits_{\R}\left(\varphi(x+y)-\varphi(x)-
y\nabla \varphi(x)\;\indyk{|y|<1}
\right)\,
  \nu(dy)\nonumber \\
&=& \lim_{\varepsilon\rightarrow 0^+} 
\int\limits_{|y|>\varepsilon}
\left(\varphi(x+y)-\varphi(x)\right)\,
  \nu(dy)\nonumber\,.
\end{eqnarray}
$\gener$ is a restriction of the infinitesimal generator of
$\{P_t\}$ \cite[Example 4.1.12]{Jc1}, and what we refer to as the 
anisotropic fractional Laplacian in the title of the paper. 
In this connection we recall that in  the special case of
$\nu(dy)=c|y|^{-d-\alpha}dy$ one obtains the fractional Laplacian
$\Delta^{\alpha/2}$. For properties of $\Delta^{\alpha/2}$ and a discussion of
equivalent definitions of its harmonic functions 
we refer the reader to \cite{BBsm1999}.

Harmonic functions corresponding to $\gener$, or $\nu$, are defined by
the mean value property with respect to an appropriate family of harmonic measures,
see Section~\ref{s:p}. 
The main goal of the paper is to characterize those operators $\gener$
for which Harnack's inequality holds, i.e., there is a constant
$C=C(\alpha,\nu)$ such that for every
function $u$ which is harmonic in the unit ball and nonnegative 
in $\R$
\begin{equation}\label{nHarnacka}
    u(x_1)\leq C u(x_2)\,,\quad |x_1|<1/2\,,\;|x_2|<1/2\,.
\end{equation}
To this end we use the relative Kato condition (RK) meaning
that there is a constant $K$ such that
\begin{equation}\label{RK}
  \int\limits_{B(y,1/2)} |y-v|^{\alpha-d} \nu(dv)
  \leq  K \nu (B(y,1/2))\,,\qquad y\in\R\,.
\end{equation}
\begin{twierdzenie}\label{Harnack}
Harnack's inequality holds for $ \gener$ if and only if {\rm (RK)} holds for $\nu$.
\end{twierdzenie}
Theorem~\ref{Harnack} is a strengthening of \cite[Theorem~1]{BS},
where an additional technical assumption was made: $\nu(dy)\leq
c|y|^{-d-\alpha}dy$, to guarantee the boundedness of $V$ on $\sfera$. 
We now drop the assumption and the boundedness is obtained 
as the sole consequence of (\ref{RK}) via Theorem~\ref{t:V}. 
We also adapt some of our previous techniques from \cite{BS} to handle measures $\nu$
which are not absolutely continuous with respect to the Lebesgue
measure on $\R$ (see, e.g., (\ref{e:P})).

Our estimates of the semigroup in Section~\ref{section:EG} are
based in part on 
ideas of \cite{Pi}, which concerns more
complicated non-convolutional semigroups.
Another, recent paper \cite{W} gives involved estimates 
of our convolution semigroup $\{P_t\}$ in individual directions 
(see also \cite{Hi} in this connection). Here we only
need isotropic estimates of $\{P_t\}$ from above, 
and our considerations become simpler than those of 
\cite{W} and \cite{Pi}.

In Section~\ref{s:p}, \ref{s:n} and \ref{s:s}
we develop the methods of \cite{BS}.
That (\ref{RK}) implies (\ref{nHarnacka}) is proved 
by using a maximum principle for a Dynkin-type version of the operator
$\gener$ to explicitly estimate its Green function $G(x,v)$ for the
unit ball, see Proposition~\ref{OG} below.
Noteworthy, our proof of the estimate is exclusive to non-local operators, of which 
$\gener$ is an illustrative special case.
In particular it turns out that $G(x,v)$ has the singularity at the
pole comparable to that of the Riesz
kernel: $|v-x|^{\alpha-d}$. 
The singularity influences the magnitude of the corresponding Poisson kernel
of the ball, $P(x,y)$, as given by the Ikeda-Watanabe formula
(\ref{e:P}). The influence is  critical if and only if (\ref{RK}) fails
to hold. This relates (\ref{RK}) to (\ref{nHarnacka}).
Such a direct influence of the singularity of the potential kernel on
the Poisson kernel does not occur for second order elliptic operators, 
which is why we can expect analogues of Theorem~\ref{Harnack} 
only for nonlocal operators.

The recent development in the study of Harnack's inequality for
general integro-differential operators similar to $\gener$ was 
initiated in \cite{BsLn}, see also \cite{BSS1}. 
The class of considered operators gradually extended, see
\cite{SV}, \cite{RSV}, \cite{BsKn},  \cite{BS}, \cite{BrBsKs}, and the references given
there. 
We note that the operators dealt with in these papers are not translation invariant 
nor are they homogeneous. On the other hand 
the papers focus on sufficient conditions 
for Harnack's inequality and they are restricted by certain 
isotropic estimates of the operator's kernel from below.

Our confinement to translation invariant homogeneous operators $\gener$ results in part
from the fact that the problem of the construction of the semigroup
from a general nonlocal operator satisfying the positive maximum principle does not
have a final solution yet. We refer the reader to \cite{S1, S2},
\cite{Jc1, Jc2}, \cite{BsKn}, and \cite{Ho}.
A general survey of the subject and more references can be found in
\cite{Bs2, JS, Jc2}.
We refer the reader to \cite{GT, Bs1} for an account of the related potential 
theory of second order elliptic operators.
We like to point out that while a symmetric second order elliptic
operator with constant coefficients is merely a linear transformation
of the Laplacian, the operators $\gener$
and their harmonic functions considered here 
are very diverse (\cite{BS}).

The remainder of the paper is organized as follows.  First definitions are
given in Section~\ref{Prel}.  In Section~\ref{section:EG} we estimate
the semigroup (see (\ref{e:p_1}) below) and the potential measure $\V$ 
and we prove our first two theorems. 
In Section~\ref{s:p} we give preliminaries needed for the proof of Theorem~\ref{Harnack},
which is presented in Section~\ref{s:n} and \ref{s:s}. In
Section~\ref{s:s} we also recall after \cite{BS} two explicit examples 
to show how irregular the L\'evy
measure $\nu$ can be for Harnack's inequality to hold or to fail for
$\gener$.

At the end of the paper we mention some remaining open problems.


\section{Preliminaries}\label{Prel}

For $x\in\R$ and $r>0$ we let $|x|=\sqrt{\sum_{i=1}^d x_i^2}$ and 
$B(x,r)=\{ y\in\R :\: |y-x|<r \}$. We denote $\sfera=\{ x\in\R:\: |x|=1 \}$.
All the sets, functions and measures considered in the sequel will be Borel.
For a measure $\lambda$ on $\R$, $|\lambda|$ denotes
its total mass.  For a function $f$ we let $\lambda(f)=\int f d\lambda$,
whenever the integral makes sense.
When $|\lambda|<\infty$ and $n=1,2,\ldots$ we let $\lambda^n$
denote the $n$-fold convolution of $\lambda$ with itself:
$$
  \lambda^n(f)=\int f(x_1+x_2+\dots+x_n)\lambda(dx_1)\lambda(dx_2)\ldots\lambda(dx_n)\,.
$$
We also let $\lambda^0=\delta_0$, the evaluation at $0$.
We call $\lambda$ degenerate if there is a proper
linear subspace $M$ of $\R$ such that $\supp(\lambda)\subset M$; 
otherwise we call $\lambda$ {\it nondegenerate}.

In what follows we will consider measures $\mu$ concentrated on $\sfera$.
We will assume that $\mu$ is positive,  finite, nondegenerate (in
particular $\mu\neq 0$), and symmetric:
$$
\mu(D)=\mu(-D)\,,\quad D\subset \R\,.
$$
We will call $\mu$ the {\it spectral measure}.
We let
\begin{equation}\label{e:mn}
  \nu(D) = \int_{\sfera} \int_0^\infty \indyk{D}(r\xi) r^{-1-\alpha} \,dr \mu(d\xi)\,,\quad
  D\subset \R\,,
\end{equation}
where $\indyk{D}$ is the indicator function of $D$. 
Note that $\nu$ is symmetric. It is a L\'evy measure on $\R$, i.e.
$$
\int_{\R}\min(|y|^2,\,1)\,\nu(dy)<\infty\,.
$$

For $r>0$ and a function $\varphi$ on $\R$ we consider its dilation
$\varphi_r(y) = \varphi(y/r)$, and we note that $\nu(\varphi_r)=r^{-\alpha} \nu(\varphi)$.
In particular $\nu$ is homogeneous: $\nu(rB)=r^{-\alpha}\nu(B)$ for $B\subset \R$.
Similarly, if $\varphi\in C^\infty_c(\R)$, 
then $\gener(\varphi_r)=r^{-\alpha}(\gener \varphi)_r$. 
This is the {\it homogeneity} of $\gener$.
In connection with the rest of our statement in Abstract we recall that
every operator $A$ on $C^\infty_c(\R)$, which satisfies the positive maximum principle:
$$
\sup_{y\in \R} \varphi(y)=\varphi(x)\geq 0 \quad  implies  \quad A \varphi(x)\leq 0\,,
$$
is given uniquely in the form
\begin{eqnarray*}
A\varphi(x)&=&\sum_{i,j=1}^d a_{ij}(x)D_{x_i}D_{x_j} \varphi(x)
+b(x) \nabla \varphi(x) -c(x)\varphi(x)\\
&&+ 
\int\limits_{\R}\left(\varphi(x+y)-\varphi(x)-
y\nabla \varphi(x)\;\indyk{|y|<1}
\right)\,
  \nu(x,dy)\,.
\end{eqnarray*}
Here $y\nabla \varphi$ is the scalar product of $y$ and the gradient
of $\varphi$
and, for every $x$, $a(x)=(a_{ij}(x))_{i,j=1}^n$ is a nonnegative definite real symmetric matrix,
the vector $b(x)=(b_i(x))_{i=1}^d$ has real coordinates, $c(x)\geq 0$,
and $\nu(x,\cdot)$ is a L\'evy measure. 
This description is due to Courr\'ege, see \cite[Proposition 2.10]{Ho},
\cite[Chapter 2]{S2} or  \cite[Chapter 4.5]{Jc1}.
For translation invariant operators $A$ the characteristics $a$, $b$, $c$,
and $\nu$ are independent of $x$. If $A$ is symmetric:
$$
\int_{\R} A\varphi(x) \phi(x)\, dx
=
\int_{\R} A\phi(x) \varphi(x)\, dx\quad \mbox{\it for} \quad \varphi,\,
\phi \in C^\infty_c(\R)\,,
$$
then $b=0$ and $\nu$ is necessarily symmetric (see, e.g., \cite[p. 251]{Jc1}
and \cite[Corollary 2.14]{Ho}).
If $A$ is homogeneous but not a local operator (\cite{Jc1}) then $a=0$
and $\nu$ must be
homogeneous, hence (\ref{e:mn}) holds with some
$\alpha\in (0,2)$ (note that $\gener \varphi(0)= \nu(\varphi)$ if 
$\varphi\in C^\infty_c(\R\setminus \{0\})$). 

We now construct the corresponding semigroup (for a more axiomatic
introduction to convolution semigroups we refer the reader to \cite{BgFt, Jc1}).
For $\varepsilon>0$ we let $\hat{\nu}_{\varepsilon}=\indyk{B(0,\varepsilon)^c}\nu$, i.e.
$\hat{\nu}_{\varepsilon}(f)=\nu(\indyk{B(0,\varepsilon)^c}f)$, and we let
$\tilde{\nu}_\varepsilon=\indyk{B(0,\varepsilon)}\nu$.
We consider the probability
measures
\begin{eqnarray}\label{eexp}
  \hat{P}^\varepsilon_t 
  &=& \exp(t(\hat{\nu}_\varepsilon-|\hat{\nu}_\varepsilon|\delta_0))
      =  \sum_{n=0}^\infty
\frac{t^n\left(\hat{\nu}_\varepsilon-|\hat{\nu}_\varepsilon|\delta_0)\right)^n}{n!}\\
&=& e^{-t|\hat{\nu}_\varepsilon|}
    \sum_{n=0}^\infty
    \frac{t^n\hat{\nu}_\varepsilon^n}{n!}\,,\quad t>0\,. \nonumber
\end{eqnarray}
Here $\hat{\nu}_\varepsilon^n=(\hat{\nu}_\varepsilon)^n$.
$\hat{P}^\varepsilon_t$ form a convolution semigroup:
$$
\hat{P}^\varepsilon_t*\hat{P}^\varepsilon_s=\hat{P}^\varepsilon_{s+t}\,,\quad s,\,t> 0\,.
$$
The Fourier transform of $\hat{P}^\varepsilon_t$ is
$$
  \Fourier(\hat{P}^\varepsilon_t)(u)=
  \int e^{iuy} \hat{P}^\varepsilon_t (dy) =
  \exp \left(t\int (e^{iuy}-1)\hat{\nu}_\varepsilon(dy)\right),
  \quad u\in\R\,.
$$
The measures $\hat{P}^\varepsilon_t$ weakly
converge to a probability measure $P_t$ as $\varepsilon\to 0$ 
(this essentially depends on (\ref{e:w}) below).
$\{ P_t, t> 0 \}$ is also a convolution semigroup and
$\Fourier(P_t)(u)=\exp(-t\Phi(u))$, where
\begin{eqnarray*}
  \Phi(u)
  &  =  & -\int \left(e^{iuy}-1-iuy\indyk{B(0,1)}(y)\right)\nu(dy) \\
   &  =  & -\int \left(\cos(uy)-1\right)\nu(dy)
=\frac{\pi}{2\sin\frac{\pi\alpha}{2}\Gamma(1+\alpha)}
  \int_{\sfera}|u\xi|^\alpha\mu(d\xi)\,.
\end{eqnarray*}
Since $\mu$ is finite and nondegenerate,
\begin{equation}\label{Phi}
\Phi(u)=|u|^\alpha \Phi(u/|u|)\approx |u|^{\alpha}\,.
\end{equation}
We call $\nu$ the {\it L\' evy measure} of the semigroup $\{P_t\,,\;t\geq
0\}$ \cite{Ho, BgFt}.

By a similar limiting procedure we construct the semigroup
$\{\tilde{P}^{\varepsilon}_t,\; t>0\}$
such that
\begin{displaymath}
  \Fourier(\tilde{P}^{\varepsilon}_t)(u) =
  \exp\left(t \int (e^{iuy}-1-iuy\indyk{B(0,1)}(y))
  \tilde{\nu}_\varepsilon(dy) \right)\,.
\end{displaymath}
Note that
\begin{equation}\label{e:w}
 \int_{\R} |y|^2 \tilde{P}_t^\varepsilon(dy)=t\int_{\R} |y|^2 \tilde{\nu}_\varepsilon(dy)\,.
\end{equation}
The L\'evy measures of   $\{\tilde{P}^\varepsilon_t\}$ and
$\{\hat{P}^\varepsilon_t\}$ are $\tilde{\nu}_\varepsilon$
and $\hat{\nu}_\varepsilon$, respectively, and we have
\begin{equation}\label{erozklad}
P_t=\tilde{P}^\varepsilon_t \ast \hat{P}^\varepsilon_t\,.
\end{equation}
The measures $P_t$ and $\tilde{P}^{\varepsilon}_t$
have rapidly decreasing Fourier transform hence they are absolutely continuous with
bounded smooth densities denoted $p_t(x)$
and $\tilde{p}^\varepsilon_t(x)$, respectively.
Of course,
\begin{equation}\label{e:wpt}
p_t=
\tilde{p}^\varepsilon_t
*
\hat{P}^\varepsilon_t
\,.
\end{equation}
By using (\ref{Phi}) we obtain the {\it scaling property} of $\{p_t\}$:
\begin{equation}\label{e:sc}
p_t(x)=t^{-d/\alpha}p_1(t^{-1/\alpha}x)\,,\quad x\in \R\,.
\end{equation}
In particular,
\begin{equation}\label{ogr_p}
  p_t(x)\leq c t^{-d/\alpha}\,.
\end{equation}
We define the {\it potential measure} of the semigroup $\{ P_t\}$:
$$
  \V(D)=\int_0^\infty P_t(D) dt\,,\quad D\subset \R\,.
$$
By (\ref{ogr_p}), $\V$ is finite on bounded subsets of $\R$ 
if $d>\alpha$. 
Let
\begin{equation}\label{eq:pkkk}
  V(x)=\int_0^\infty p_t(x)dt\,,\quad x\in\R\,,
\end{equation}
so that
$$
 \V(D)=\int_D V(x)dx\,,\;\quad D\subset \R\,.
$$
We call $V(x)$ the {\it potential kernel} of the stable semigroup.
By (\ref{e:sc}) 
\begin{equation}
  \label{e:sv}
  V(x)=|x|^{\alpha-d} V\left(x/|x|\right)\,,\quad x\neq 0\,,
\end{equation}
and $\V(rD)=r^\alpha \V(D)$ for $r>0$, $D\subset\R$.

If $d=1$ then up to a constant there is only one measure $\nu$ to
consider: $\nu(dy)=|y|^{-1-\alpha}dy$, corresponding to 
 $\gener=c\Delta^{\alpha/2}$. 
This case is not excluded from our considerations but it is sometimes trivial.
In particular, if $d=1\leq \alpha$ then  $V\equiv\infty$
(\cite[Example 14.30]{BgFt}).
We refer to \cite{BZ} for more information and references on the case $d=1\leq \alpha$.

{\it Constants} in this paper mean positive real numbers.  
We often write $f\approx g$ to indicate that there is  
$c=c(\alpha,\mu)$, i.e. a constant $c$ depending only on $\alpha$ and $\mu$, 
such that $c^{-1}f \leq g \leq c f$.


\section{Estimates of semigroup and potential measure}
\label{section:EG}
A general reference to the potential theory of convolution semigroups is
\cite{BgFt} (see also \cite{Jc1, Jc2}).

We consider an auxiliary scale of smoothness for $\nu$.
\begin{definicja}
We say that $\nu$ is a $\gamma$-measure on $\sfera$ if
\begin{equation}\label{nu_gamma} 
  \nu(B(x,r)) \leq c r^{\gamma}\,,\quad |x|=1\,,\; 0<r<1/2\,.
\end{equation}
\end{definicja}
Since $\nu(drd\theta)=r^{-1-\alpha}dr\mu(d\theta)$, it is at least a $1$-measure
and at most a $d$-measure on $\sfera$.
If $\nu$ is a $\gamma$-measure with $\gamma>1$, then $\mu$ has no atoms.
$\nu$ is a $d$-measure if and only if it is absolutely continuous with
respect to the Lebesgue measure and has a density function which is locally bounded
on $\R\setminus\{0\}$. 
We refer the reader to \cite{Dz} and \cite{GH} for considerations related to this case.

In the remainder of this section we fix $1\leq \gamma \leq d$ and we assume
that $\nu$ is a $\gamma$-measure on $\sfera$.

\noindent
We will first estimate individual terms in the series in (\ref{eexp}).
\begin{lemat}\label{Conv_est}
There exists
$C=C(\alpha,\mu)$ such that for $\varepsilon>0$ and $n=1,2,\ldots$  we have
\begin{equation}\label{e:cest}
  \nusmall{\varepsilon}^n(B(x,r)) \leq
  C^n r^\gamma \varepsilon^{-(n-1)\alpha}\,,\quad |x|=1\,,
\end{equation}
provided $0<r<\max(\varepsilon/3,1/5^n)$.
\end{lemat}
\dowod 
We proceed by induction. Note that (\ref{e:cest}) holds for $n=1$ by (\ref{nu_gamma}).
Let $c_0$ and $n$ be such that (\ref{e:cest}) is satisfied with $C=c_0$.
We first assume that $r<\varepsilon/3$. For every $x\in\sfera$ by homogeneity of $\nu$ and (\ref{nu_gamma}) we have
\begin{eqnarray*}
  \nusmall{\varepsilon}^{n+1}(B(x,r))
  &   =  & \int\limits_{|x-y|>2\varepsilon/3} \nusmall{\varepsilon}(B(x-y,r))
           \nusmall{\varepsilon}^n(dy) \\
  & \leq & \int\limits_{|x-y|>2\varepsilon/3} \nu(B(x-y,r))\nusmall{\varepsilon}^n(dy) \\
  &   =  & \int\limits_{|x-y|>2\varepsilon/3}
           |x-y|^{-\alpha}
           \nu(B(\frac{x-y}{|x-y|},\frac{r}{|x-y|}))\nusmall{\varepsilon}^n(dy) \\
  & \leq & c_1 r^{\gamma} \int\limits_{|x-y|>2\varepsilon/3}
            |x-y|^{-\alpha-\gamma} \nusmall{\varepsilon}^n(dy)
\end{eqnarray*}
(note that $r/|x-y|<1/2$ provided $|x-y|>2\varepsilon/3$).
Now let $\varepsilon/3\leq r<1/5^{n+1}$. Then $2r+\varepsilon<1/5^n$ and 
by induction
\begin{eqnarray*}
  \int\limits_{|x-y|<2r+\varepsilon}\nusmall{\varepsilon}(B(x-y,r))
  \nusmall{\varepsilon}^n(dy)
  & \leq & |\nusmall{\varepsilon}| \nusmall{\varepsilon}^n(B(x,2r+\varepsilon)) \\
  & \leq & \frac{|\mu|}{\alpha} \varepsilon^{-\alpha}
           c_0^{n}(2r+\varepsilon)^\gamma \varepsilon^{-(n-1)\alpha} \\
  & \leq & c_0^n c_2 r^{\gamma} \varepsilon^{-n\alpha},
\end{eqnarray*}
for some $c_2=c_2(\alpha,\mu)$ ; and by homogeneity of $\nu$ and (\ref{nu_gamma}) we get
\begin{eqnarray*}
  \int\limits_{|x-y|>2r+\varepsilon}\nusmall{\varepsilon}(B(x-y,r))
  \nusmall{\varepsilon}^n(dy)
  & \leq & \int\limits_{|x-y|>2r+\varepsilon}\nu(B(x-y,r))
           \nusmall{\varepsilon}^n(dy) \\
  & \leq & \int\limits_{|x-y|>2r+\varepsilon} c_1 r^{\gamma}
           |x-y|^{-\alpha-\gamma} \nusmall{\varepsilon}^n(dy) \\
  & \leq & c_1 r^{\gamma}\int\limits_{|x-y|>2\varepsilon/3} 
           |x-y|^{-\alpha-\gamma} \nusmall{\varepsilon}^n(dy).
\end{eqnarray*}
From the above we have
\begin{equation}\label{podsum}
  \nusmall{\varepsilon}^{n+1}(B(x,r)) \leq c_1 r^{\gamma}\int\limits_{|x-y|>2\varepsilon/3} 
           |x-y|^{-\alpha-\gamma} \nusmall{\varepsilon}^n(dy)+c_0^n c_2 r^{\gamma} \varepsilon^{-n\alpha},
\end{equation}
for all $0<r<\max(\varepsilon/3,1/5^{n+1})$.

Let $L_\varepsilon=\left\lfloor \log_5(3/2\varepsilon)\right\rfloor$. If $2\varepsilon/3<1/5^n$ then we get by induction
$$
\begin{array}{lcl}
  \int\limits_{2\varepsilon/3<|x-y|<1/5^n}
  |x-y|^{-\alpha-\gamma} \nusmall{\varepsilon}^n(dy)
  & \leq & \sum\limits_{k=n}^{L_\varepsilon}\; \int\limits_{1/5^{k+1}<|x-y|<1/5^k}
           |x-y|^{-\alpha-\gamma} \nusmall{\varepsilon}^n(dy) \\
  & \leq & \sum\limits_{k=n}^{L_\varepsilon}
           (5^{k+1})^{\alpha+\gamma} \nusmall{\varepsilon}^n(B(x,1/5^k)) \\
  & \leq & c_0^n 5^{\alpha+\gamma} \varepsilon^{-(n-1)\alpha}
           \sum\limits_{k=1}^{L_\varepsilon} 5^{k\alpha} \\
  & \leq & c_0^n c_3 \varepsilon^{-n\alpha}\,,
\end{array}
$$
where $c_3=c_3(\alpha,\mu)$.
Also,
$$
\begin{array}{lcl}
  \int\limits_{|x-y|>1/5^n} |x-y|^{-\alpha-\gamma}\nusmall{\varepsilon}^n(dy)
  & \leq & (5^{\alpha+\gamma})^n |\nusmall{\varepsilon}^n| \\
  &   =  & (5^{\alpha+\gamma} \frac{|\mu |}{\alpha})^n \varepsilon^{-n\alpha} \\
  & \leq & c_0^n \varepsilon^{-n\alpha},
\end{array}
$$
by taking large $c_0$.
We get
$$
  \int\limits_{|x-y|>2\varepsilon/3}
  |x-y|^{-\alpha-\gamma} \nusmall{\varepsilon}^n(dy) \leq c_0^n \varepsilon^{-n\alpha}(c_3+1),
$$
and (\ref{podsum}) yields
$$
  \nusmall{\varepsilon}^{n+1}(B(x,r)) \leq c_0^{n+1} r^{\gamma}\varepsilon^{-n\alpha}.
$$ 
\qed

\begin{wniosek}\label{lambda_epsilon}
There exists
$C=C(\alpha,\mu)$ such that
\begin{equation}\label{eq_lambda_epsilon}
  \nusmall{\varepsilon}^n(B(x,\lambda\varepsilon)) \leq
  C^n \lambda^{\gamma} (1+\lambda^{\alpha})
  \varepsilon^{\gamma-(n-1)\alpha}\,,\quad \lambda>0\,,\; \varepsilon>0\,,\; |x|=1\,.
\end{equation}
\end{wniosek}
\dowod Lemma \ref{Conv_est} yields (\ref{eq_lambda_epsilon}) for
$\lambda\varepsilon<1/5^n$. For $\lambda\varepsilon\geq 1/5^n$ we have
$$
  \nusmall{\varepsilon}^n(B(x,\lambda\varepsilon)) \leq |\nusmall{\varepsilon}^n|
  = \frac{|\mu|^n}{\alpha^n}\varepsilon^{-n\alpha} \leq
  (\frac{|\mu|}{\alpha}5^{\alpha+\gamma})^n\lambda^{\alpha+\gamma}
  \varepsilon^{\gamma-(n-1)\alpha}.
$$
\qed

\noindent
In what follows we denote $\hat{P}_t=\hat{P}_t^{t^{1/\alpha}}$
and $\tilde{P}_t=\tilde{P}_t^{t^{1/\alpha}}$.

\begin{wniosek}\label{Sem_est}
There exists
$C=C(\alpha,\mu)$ such that
$$
  \hat{P}_t(B(x,\lambda t^{1/\alpha})) \leq
  C \lambda^{\gamma}(1+\lambda^{\alpha})
  t^{1+\frac{\gamma}{\alpha}}\,,\quad  \lambda>0\,,\; t>0\,,\; |x|=1.
$$
\end{wniosek}
\dowod Corollary \ref{lambda_epsilon} yields
\begin{eqnarray*}
  \hat{P}_t(B(x,\lambda t^{1/\alpha}))
  &  =   & e^{-|\mu|/\alpha} \sum\limits_{n=0}^\infty
           \frac{t^n \nusmall{t^{1/\alpha}}^n(B(x,\lambda t^{1/\alpha}))}{n!} \\
  & \leq & e^{-|\mu|/\alpha}
           \sum\limits_{n=0}^\infty
           \frac{c^n
           \lambda^\gamma(1+\lambda^\alpha)t^{1+\frac{\gamma}{\alpha}}}
           {n!} \\
  &   =  & e^{c-|\mu|/\alpha}
           \lambda^\gamma(1+\lambda^\alpha)t^{1+\frac{\gamma}{\alpha}}.
\end{eqnarray*}
\qed

\begin{wniosek}\label{w:p1h}
$\hat{P}_1(B(y,\lambda))\leq C
\lambda^\gamma(1+\lambda^\alpha)|y|^{-\alpha-\gamma}$ for
$y\in \R$ and $\lambda>0$.
\end{wniosek}
\dowod
Let $y\in \R\setminus\{0\}$ and $x=y/|y|$, $t=|y|^{-\alpha}$.
By scaling and Corollary~\ref{Sem_est} we have\\
$\hat{P}_1(B(y,\lambda))=
\hat{P}_t(B(x,\lambda t^{1/\alpha}))
\leq c \lambda^\gamma(1+\lambda^\alpha)|y|^{-\alpha-\gamma}\,.
$
\qed

\vspace{10pt}
We note that for every $q>0$ we have that $\int |y|^q \tilde{P}_1(dy)<\infty$,
because the support of $\tilde{\nu}_1$ is bounded (\cite{Sato}).
A simple reasoning  based on this and the boundedness of
the derivative of $\tilde{p}_1$ 
yields
$$
\tilde{p}_1(y)\leq c_q (1+|y|)^{-q}\,,\quad q>0\,,\; y\in\R\,,
$$
see \cite[Lemma 9]{Pi}.
\begin{lemat}\label{l:p1f}
For every $q>0$ there exists $C=C(\alpha,\mu,q)$ such that
$$
\tilde{P}_1(B(z,\rho))\leq C
(1+|z|)^{-q}\rho^d\,, \quad \rho\leq
1\,,\; z\in \R\,.
$$
\end{lemat}
\dowod
If $|z|<2$ then
$\tilde{P}_1(B(z,\rho)) =\int_{B(z,\rho)} \tilde{p}_1(y)dy
\leq c \rho^d \leq c (1+|z|)^{-q}\rho^d$.\\
If $|z|\geq 2$ then
$\tilde{P}_1(B(z,\rho)) \leq c (1+|z|/2)^{-q}\rho^d \leq c(1+|z|)^{-q}\rho^d$.
\qed

The proof of the following lemma is a simplification of the proof of \cite[Theorem 3]{Pi}.
\begin{lemat}\label{l:p1}
$P_1(B(z,\rho))\leq C|z|^{-\alpha-\gamma}\rho^d$ for $z\in \R$ and $0<\rho\leq 1$.
\end{lemat}
\dowod
By (\ref{erozklad}), Lemma~\ref{l:p1f}, and Corollary~\ref{w:p1h}
\begin{eqnarray*}
  P_1(B(z,\rho))
  &  =   & \tilde{P}_1*\hat{P}_1(B(z,\rho)) = \int_{\R}\tilde{P}_1(B(z-y,\rho))\hat{P}_1(dy) \\
  &  =   & \int_0^1 \hat{P}_1(\{y\,:\; \tilde{P}_1(B(z-y,\rho))>s\})ds \\
  & \leq & \int_0^1 \hat{P}_1(\{y\,:\; c(1+|z-y|)^{-q} \rho^d>s\})ds \\
  & \leq & \int_0^{c\rho^d} \hat{P}_1(B(z, c^{1/q}s^{-1/q}\rho^{d/q}))ds\\
  & \leq & c \int_0^{c\rho^d} (c^{1/q}s^{-1/q}\rho^{d/q})^\gamma
             (1+(c^{1/q}s^{-1/q}\rho^{d/q})^\alpha)|z|^{-\gamma-\alpha}ds \\
  &  =   & c |z|^{-\gamma-\alpha} \left[
              \rho^{d\gamma/q}\int_0^{c\rho^d} s^{-\gamma/q}ds
              +
              \rho^{d(\gamma+\alpha)/q}\int_0^{c\rho^d} s^{-(\gamma+\alpha)/q}ds
             \right] \\
  &  =   & c |z|^{-\gamma-\alpha} \left[
              \rho^{d\gamma/q} (\rho^d)^{1-\gamma/q}
              +
              \rho^{d(\gamma+\alpha)/q}(\rho^d)^{1-(\gamma+\alpha)/q}
              \right]
             = c |z|^{-\gamma-\alpha}\rho^d\,.
\end{eqnarray*}
\qed

The following two corollaries are our main estimates of the semigroup.
Corollary~\ref{w:ptrho} is an analogue of \cite[Theorem 3]{Pi}, 
while (\ref{e:p_1}) corresponds to \cite{W}.
\begin{wniosek}\label{w:rho}
$P_1(B(z,\rho))\leq C(1+|z|)^{-\alpha-\gamma}\rho^d\;$ if
$\;0\leq\rho< |z|/2$.
\end{wniosek}
\dowod
We recall that $p_1(y)=P_1(dy)/dy$ is bounded and so Lemma~\ref{l:p1} yields 
\begin{equation}\label{e:p_1}
  p_1(y)\leq c(1+|y|)^{-\gamma-\alpha}\,,\quad y\in \R\,.
\end{equation}
If $0\leq \rho<|z|/2$
then
$
P_1(B(z,\rho))\leq c\int_{B(z,\rho)}(1+|y|)^{-\gamma-\alpha}dy
\leq(1+|z|)^{-\alpha-\gamma}\rho^d\,.
$
\qed

\begin{wniosek}\label{w:ptrho}
$
P_t(B(x,\rho))\leq C t^{1+\frac{\gamma-d}{\alpha}}\rho^d
$,
provided $|x|=1$, $t>0$, and $0\leq \rho\leq t^{1/\alpha}$.
\end{wniosek}
\dowod
By scaling and Lemma~\ref{l:p1} we have
$
P_t(B(x,\rho))=
P_1(B(xt^{-1/\alpha},\rho t^{-1/\alpha}))
\leq c t^{1+\frac{\gamma-d}{\alpha}}\rho^d
$.
\qed

\vspace{5mm}
\noindent
{\it Proof of {\rm Theorem~\ref{t:V}}}. 
Let $|x|=1$, $0\leq \rho<1/2$. By scaling and Corollary~\ref{w:rho}
\begin{eqnarray*}
  \V(B(x,\rho))&=&\int_0^\infty P_t(B(x,\rho))dt
=\int_0^\infty P_1(B(xt^{-1/\alpha},\rho t^{-1/\alpha}))dt\\
&\leq&
c\rho^d \int_0^\infty (1+t^{-1/\alpha})^{-\gamma-\alpha}
t^{-d/\alpha}dt\,.
\end{eqnarray*}
The integral is finite because $-d/\alpha<-1$ and $(\gamma+\alpha-d)/\alpha>-1$.
Let $y\in \R\setminus \{0\}$, $x=y/|y|$. By
scaling, a change of variable, and (\ref{e:p_1})
\begin{eqnarray*}
V(y)&=&\int_0^\infty t^{-d/\alpha}p_1(yt^{-1/\alpha})dt
=|y|^{\alpha-d} \int_0^\infty t^{-d/\alpha}p_1(xt^{-1/\alpha})dt \nonumber\\
&\leq &
|y|^{\alpha-d} \int_0^\infty t^{-d/\alpha}(1+t^{-1/\alpha})^{-\gamma-\alpha}dt \leq
c|y|^{\alpha-d}\,.
\label{e:v}
\end{eqnarray*}
The first integral above is locally uniformly convergent on
$\R\setminus\{0\}$ hence $V$ is continuous there.
\qed

We now proceed to our converse, Theorem~\ref{t:nu}.
We propose a general approach
based on a simple study of generator $\gener$.
We first note that 
\begin{equation}\label{typeA}
  p_t(x)>0\,, \quad x\in\R\; \quad (t>0)\,,
\end{equation}
see (\cite{Taylor}) or \cite[Lemma 5]{Pi}. 
In fact, (\ref{typeA}) easily follows from (\ref{e:wpt}), (\ref{eexp}), continuity
of $\tilde{p}^\varepsilon_t$, and the fact that 
$\supp(\nu)+\supp(\nu)+\ldots+\supp(\nu)$ ($d$ times) equals $\R$.

By (\ref{typeA}), (\ref{e:sv}), and continuity of $p_t$ for $t>0$, there is
a constant $c=c(\alpha,\mu)$ such that
\begin{equation}\label{v>C}
V(x) \geq c |x|^{\alpha-d},\; x\in \R.
\end{equation}

\begin{lemat}\label{lemat:generator2}
Let $d>\alpha$. For all $\varphi\in C_c^\infty(\R)$ we have
$$
  \int\limits_{\R} \gener\varphi(x-y)\V(dy)=-\varphi(x)\,,\quad x\in\R\,,
$$
where the integral is absolutely convergent.
\end{lemat}
This is well-known (see, e.g., \cite[Theorem 3.5.78]{Jc2}). We only 
note that $|\gener \varphi(x)|\leq c(1+|x|)^{-1-\alpha}$.
The absolute convergence follows from this and the homogeneity of $\V$.

\vspace{5mm}

\noindent
{\it Proof of {\rm Theorem~\ref{t:nu}}.}
If $d=1\leq \alpha$ then $V\equiv\infty$ and there is nothing to prove.
Thus we assume that $d>\alpha$.
We fix a function $\phi\in C_c^\infty(\R)$ such that  $\phi\geq
0$, ${\rm supp}\, \phi\subset B(0,1/2)$ and $\phi=1$ on $B(0,1/3)$.
Let $r>0$.  Put $\phi_r(x)=\phi(x/r)$ and
$\Lambda_r(x)=\gener \phi_r(x)$.
Homogeneity of $\gener$  yields
$
  \Lambda_r(x)=r^{-\alpha}\Lambda_1(x/r)
$.
Note that $\gener \phi=\Lambda_1$ is bounded, hence
there is a constant $c$ such that
$$  \Lambda_r(x)\geq-cr^{-\alpha}\,.
$$
If $|x|\geq r/2$ then $\Lambda_r(x)\geq 0$, and in fact
$\Lambda_r(x)\geq \nu(B(x,r/3))$.
Let $|x|>r$. From Lemma~\ref{lemat:generator2} we have
\begin{eqnarray*}
  0&=&
 \int\limits_{\R} \Lambda_r(x-y)\V(dy)
\geq
 \int\limits_{B(x,r/2)} \Lambda_r(x-y)\V(dy)
+
 \int\limits_{B(0,r/4)} \Lambda_r(x-y)\V(dy)\\
&\geq&
 -cr^{-\alpha}\V((B(x,r/2))
+
 \int\limits_{B(0,r/4)} \nu(B(x-y,r/3))\V(dy)\\
&\geq&
 -cr^{-\alpha}\V((B(x,r/2))
+
\V(B(0,r/4)) \nu(B(x,r/12))\,.
\end{eqnarray*}
Since $\V(B(0,r/4))=r^{\alpha}\V(B(0,1/4))$ and $\V(B(0,1/4))<\infty$ we get
\begin{equation}
  \label{e:nuV}
\nu(B(x,r/12))\leq cr^{-2\alpha}\V((B(x,r/2))\,,
\quad |x|>r\,.
\end{equation}
\qed

We note that similar results can also be derived from the lower bounds for the semigroup 
as given in \cite[Theorem 1.1]{W}. 


\section{Harnack's  inequality: preliminaries}\label{s:p}
The general references for this section are \cite{Ch, Ch1}, \cite{Sato}, \cite{Br}, or \cite{BGB}.  
The L\' evy measure $\nu$ yields a {\it standard} symmetric stable L\' evy process $\proces$
with generating triplet $(0,\nu,0)$. Namely, the transition probabilities of the process $\proces$ are 
$P(t,x,A)=P_t(A-x)$, $t>0$, $x\in\R$, $A\subset\R$, and 
$P(0,x,A)=\indyk{A}(x)$, where $\{ P_t,\; t\geq 0\}$ is the stable semigroup of 
measures introduced in Preliminaries. The process is strong Markov with respect 
to the so-called {\it standard filtration}.

The process conveniently leads to a definition of harmonic measures
$\omega^x_D$, and their properties (\ref{e:www}) and  (\ref{wzornaGreena}) below. 
For an analytic definition of these, called the fundamental
family, we refer to \cite{BgFt} (see also \cite{La, BH}).

For open $U\subset \R$ we denote $\tau_U=\inf\{t \ge 0;\: X_t\not\in U\}$, the
{\it first exit time} of $U$.
We write $\omega_D^x$ for the harmonic measure of (open) $D$:
$$
  \omega_D^x(A)=P^x(\tau_D<\infty\,,\;X_{\tau_D}\in A)\,,\quad x\in \R\,,\; A\subset \R\,.
$$
By the strong Markov property
\begin{equation}
  \label{e:www}
  \omega^x_D(A)=\int \omega^y_D(A)\omega^x_U(dy)\,,\quad\mbox{ if }\; U\subset D\,.
\end{equation}

We say that a function $u$ on $\R$ is
{\it harmonic} in  open  $D\subset \R$ if
\begin{equation}\label{wyrazenie:dH}
  u(x)=E^x u(X_{\tau_U})=\int_{U^c} u(y)\,\omega_U^x(dy)\,,\quad x\in \R,
\end{equation}
for every bounded open set $U$ with the closure $\bar{U}$ contained in
$D$.  It is called {\it regular harmonic} in $D$ if
(\ref{wyrazenie:dH}) holds for $U=D$. If $D$ is unbounded then
$E^x u(X_{\tau_D}) = E^x[\tau_D < \infty\,;\; u(X_{\tau_D})]$
by a  convention.  Under (\ref{wyrazenie:dH}) it will be only assumed that
the expectation in (\ref{wyrazenie:dH}) is well defined (but not necessarily finite).
Regular harmonicity implies harmonicity, and it is inherited by subsets $U\subset D$. 
This follows from (\ref{e:www}).

We denote by $p_t^D(x,v)$ the transition density of the process killed 
at the first exit from $D$:
$$p_t^D(x,v)=p(t,x,v)-E^x[\tau_D<t\,;\;p(t-\tau_D,X_{\tau_D},v)]\,,\quad
t>0\,,\; x,v \in \R\,.$$ 
Here $p(t,x,v)=p_t(v-x)$. 
For convenience we will assume in the sequel that $D$ is regular:
$P^x[\inf\{t>0\,:\;X_t\notin D\}=0]=1$ for $x\in D^c$ (see \cite{Ch1, Ch}).
Then $p_t^D$ is symmetric: $p_t^D(x,v)=p_t^D(v,x)$, $x, v\in D$ (see, e.g., \cite{ChZ}). 
The strong Markov property yields
\begin{equation}
  \label{e:pt}
   p(t,x,v) = E^x[p(t-\tau_D,X_{\tau_D},v)\,;\;\tau_D<t]\,,
\quad x\in D\,,\; v\in D^c\,.
\end{equation}
In particular, $p_t^D(x,v)=0$ if $x\in D$, $v\in D^c$.
We let
  \begin{displaymath}
G_D(x,v)=\int\limits_0^\infty p_t^D(x,v) dt\,,
  \end{displaymath}
and we call $G_D(x,v)$ the {\it Green function} for $D$.  
If $V$ is continuous on $\R\setminus \{0\}$, so that $V(x)\leq
c|x|^{\alpha-d}$, then the strong Markov property yields for $x,v\in D$
\begin{equation}\label{wzornaGreena}
  G_D(x,v)=
V(x,v)-E^xV(X_{\tau_D},v)
=V(x,v)-\int_{D^c}V(z,v)\,\omega_D^x(dz)\,.
\end{equation}
Here $V(x,v)=V(v-x)$.
The Green function is symmetric: $G_D(x,v)=G_D(v,x)$, continuous
in $D\times D\setminus\{(x,v): x=v\}$, and it vanishes if $x\in D^c$ or $v\in D^c$.

Note that $V(x,v)$ is harmonic in $x$ on $\R\setminus\{ v \}$. 
Indeed, if $x\in D$ and $\dist(D,v)>0$ then by (\ref{e:pt}) 
\begin{eqnarray*}
  V(x,v)
  &  =   & \int\limits_0^\infty E^x[p(t-\tau_D,X_{\tau_D},v);\tau_D<t] dt 
  =   E^x V(X_{\tau_D},v)\,.
\end{eqnarray*}
Similarly, the Green function $v\mapsto G_D(x,v)$ is harmonic in $D\setminus\{ x \}$.

By Ikeda--Watanabe formula \cite{IW} we have
\begin{equation}\label{IW}
  \omega_D^x(A) = \int_D G_D(x,v)\nu(A-v)dv\,,\quad \mbox{ if }\; \dist(A,D)>0\,.
\end{equation}
We note here that translation--invariance of the Lebesgue measure and
Fubini--Tonelli theorem yield
\begin{equation}\label{PhiPsi}
  \int \int \Phi(v) \Psi(v+z) m(dz) dv = \int \int \Phi(v+z) \Psi(v) m(dz) dv\,,
\end{equation}
for every symmetric measure $m$ and nonnegative functions $\Phi$ and $\Psi$. In particular,
taking $m=\nu$, $\Phi(v)=G_D(x,v)$ and $\Psi(v)=\indyk{A}(v)$ we get
$$
  \int_{D}G_D(x,v)\nu(A-v)dv =
  \int_A \int_{-D+v} G_D(x,v-z) \nu(dz) dv\,.
$$
If the boundary of $D$ is smooth or even Lipschitz then
$$
  \omega_{D}^x(\partial D)=0\,,\quad x\in D\,,
$$
see \cite{Sztonyk} (see also \cite{Millar}, \cite{Wu}). 
In this case $\omega_D^x$ is absolutely continuous with respect to
the Lebesgue measure on $D^c$.
Its density function, or the Poisson kernel, is given by the formula
\begin{equation}
  \label{e:P}
  P_D(x,y)=\int_{y-D} G_D(x,y-z) \nu(dz)\,,\quad x\in D\,.
\end{equation}
Note that such $D$ are regular, because of (\ref{typeA}) and scaling.
In particular the above considerations apply to $D=B(0,1)$.
 
It follows from (\ref{e:sc}) that for every $r>0$ and $x\in \R$ the
$P_x$ distribution of $\{X_t\,,\,t\geq 0\}$ is the same as the $P_{rx}$
distribution of $\{r^{-1}X_{r^{\alpha}t}\,,\,t\geq 0\}$. 
In particular,
\begin{equation}
  \label{e:scx}
\omega^x_D(A)=\omega^{rx}_{rD}(rA)\,.
\end{equation}
We call (\ref{e:scx}) scaling, too.
It yields that for $u$ harmonic on $D$, the dilation,
$u_r$, is harmonic on $rD$. A similar remark concerns translations. 
 
By (\ref{PhiPsi}) we also obtain
\begin{displaymath}
  \int_{B(0,1/2)}|y|^{\alpha-d}\nu(A-y)dy =
  \int_A \int_{B(y,1/2)} |y-z|^{\alpha-d} \nu(dz) dy\,,\quad A\subset \R\,,
\end{displaymath}
and
$$
  \int_{B(0,1/2)}\nu(A-y)dy =
  \int_A \nu(B(y,1/2)) dy\,,\quad A\subset \R\,.
$$
Therefore we can express the relative Kato condition (RK) in an equivalent form:
\begin{equation}\label{RKm}
  \int_{B(0,1/2)}|y|^{\alpha-d}\nu(A-y)dy
  \leq K \int_{B(0,1/2)}\nu(A-y)dy\,,\quad A\subset \R\,.
\end{equation}

We remark that (RK) is a local condition at infinity: the
inequality in (\ref{RK}) only needs to be verified for large $y\in \R$. 
In particular, if it holds for $|y|>1$ then it holds for all $y\in\R$,
possibly with a different constant, see \cite{BS}.
Noteworthy, the reverse of (\ref{RK}) (and (\ref{RKm})) always holds,
so actually (RK) means comparability of both sides of (\ref{RK}) (and (\ref{RKm})).

In what follows we let $G=G_{B(0,1)}$, $P=P_{B(0,1)}$ and we define
$$s(x)=E^x \tau_{B(0,1)}=\int_{B(0,1)} G(x,v)dv\,.$$
Explicit formulas for these functions for
$\nu(dy)=|y|^{-d-\alpha}dy$ are known and may give some insight into the general situation.
They are essentially due to M. Riesz, see, e.g., \cite{BB2}, \cite{BGR}, \cite{La}, \cite{BH}, \cite{G}.
In particular (for isotropic $\nu$) we have
\begin{equation}
  \label{e:jp}
P(x,y)=C^d_\alpha\left[\frac{1-|x|^2}{|y|^2-1}\right]^{\alpha/2}|x-y|^{-d}\,,
\qquad |x|<1\,, |y|> 1\,.
\end{equation}

The following two lemmas
are consequences of symmetry and nondegeneracy of the spectral measure $\mu$.
They can be proved similarly as Lemma 4 and Lemma 10 of \cite{BS}, so we skip the proofs.

\begin{lemat}\label{czescDynkina}
  There exist 
$\varepsilon = \varepsilon(\alpha,\mu)\in (0,1)$ and
  $C=C(\alpha,\mu)$ such that
\begin{equation}\label{Lduza}
  \nu (B(x,1-\varepsilon)) \geq C\,,
\end{equation}
provided $1-\varepsilon < |x| < 1$.
\end{lemat}

\begin{lemat}\label{ET_est} There exists 
$C=C(\alpha,\mu)$ such that
$$
  s(x) \leq C (1-|x|^2)^{\alpha/2}\,,\quad |x|<1.
$$
\end{lemat}

\section{Necessity of relative Kato condition}\label{s:n}
In this short section we assume that Harnack's  inequality
(\ref{nHarnacka}) holds. 
We make no further assumptions on $\nu$ beyond these in Section~\ref{Prel}.
In particular our considerations do not depend on the estimates in Section \ref{section:EG}.

\begin{lemat}\label{n_RK}
  Harnack's  inequality implies the relative Kato condition.
\end{lemat}
\dowod
We first consider the case $d>\alpha$. 
We claim that
\begin{equation}\label{e:cv}
V(x)\approx |x|^{\alpha-d}\,,\quad x\in\R.
\end{equation}
Indeed, for every $|x|=1$, 
$\V(B(x,1/4))=\int_{B(x,1/4)}V(v)dv \leq \V(B(0,2)) < \infty $,
so there exists $v\in B(x,1/4)$ such that $V(v) \leq \V(B(0,2)) /
|B(0,1/4)|$. By Harnack's inequality $V(x)\leq c V(v)$. 
The estimate (\ref{e:cv}) follows from (\ref{e:sv}) and (\ref{v>C}).

Let $g(v)=\min(G(0,v),1)$. We claim that
\begin{equation}
  \label{e:ggg}
G(x,v) \approx g(v) |v-x|^{\alpha-d} \quad if \;
\; |x|<1/2 \; and \; |v|<1\,.
\end{equation}
Indeed, by (\ref{e:cv})  and (\ref{wzornaGreena}) for small $\delta>0$
we have:
$$
  G(x,v) \approx |v-x|^{\alpha-d}\,,\quad |x|<1/2\,,\; |x-v|<\delta\,.
$$
Harnack's  inequality yields that
$ G(x,v) \approx |v-x|^{\alpha-d}$ provided $|x|<1/2$ and $|v|<3/4$, and also
$G(x,v)\approx G(0,v)$ if $|x|<1/2$ and $|v|>3/4$.
Note that $g$ is locally bounded from below on $B(0,1)$. This
completes the proof of (\ref{e:ggg}).

For every  $A\subset \R$ the function $x\mapsto\omega_{B(0,1)}^x(A)$ is nonnegative on $\R$
and regular harmonic in $B(0,1)$. 
Harnack's  inequality (\ref{nHarnacka}), (\ref{IW}), (\ref{e:ggg}), and
Fubini-Tonelli yield
\begin{eqnarray*}
  \omega_{B(0,1)}^0(A) 
  & \approx &  \int_{B(0,1/2)} \omega^x_{B(0,1)}(A)\,dx \\
  & \approx &  \int_B \int_{B(0,1/2)} g(v) |v-x|^{\alpha-d}\nu(A-v)\,dv\,dx \\
  & \approx &  \int_B g(v) \nu(A-v)\,dv\,.
\end{eqnarray*}
This and (\ref{IW}) yield
\begin{displaymath}
  \int_B g(v)|v|^{\alpha-d}\nu(A-v)\,dv \approx
  \int_B g(v)\nu(A-v)\,dv\,.
\end{displaymath}
To this ``approximate equality'' we add the following one:
\begin{displaymath}
  \int_{B\setminus B(0, 3/4)} |v|^{\alpha-d}\nu(A-v)\,dv \approx
  \int_{B\setminus B(0,3/4)} \nu(A-v)\,dv\,,
\end{displaymath}
and we obtain
\begin{displaymath}
  \int_{B} |v|^{\alpha-d}\nu(A-v)\,dv \approx
  \int_{B} \nu(A-v)\,dv\,,\quad A\subset B^c\,.
\end{displaymath}
A change of variable: $v=2u$ yields (\ref{RKm}) and (\ref{RK}). 

In the case $d\leq \alpha$ we have $d=1$, and so
$\nu(dy)=c|y|^{-1-\alpha}dy$, which satisfies (RK).
\qed

\section{Sufficiency of relative Kato condition}
\label{s:s}In what follows we assume that (RK) holds for $\nu$.
We will also assume that $d>\alpha$ unless stated otherwise.

The key step in the proof of Harnack's  inequality 
is the following estimate for the Green function of the ball, 
which we prove after a sequence of lemmas. 
We note that it is essentially the same inequality as (\ref{e:ggg}), 
but proved under explicit assumptions on $\nu$ rather than 
by stipulating Harnack's inequality.
The estimate was suggested by the sharp estimates of the Green function of Lipschitz domains \cite{Jk} 
for the isotropic $\nu$ (see also \cite{Bo2}). 
We  also refer the reader to \cite{K1, CS} for
more explicit estimates for smooth domains and to, e.g., \cite{BB2} 
for explicit formulas for the ball. 

\begin{prop}\label{OG} 
   $G(x,v) \approx s(v)|v-x|^{\alpha-d}$ provided $|x|<1/2$ and $|v|<1$.
\end{prop}

\begin{lemat}\label{nugamma} $\nu$ is a $(d-\alpha)$-measure on
  $\sfera$.
\end{lemat}
\dowod Indeed, for $|x|=1$, $0<r<1/2$
by (\ref{RK}) we obtain
$$
  \nu(B(x,r)) \leq r^{d-\alpha} \int\limits_{B(x,1/2)}|x-z|^{\alpha-d} \nu(dz)
  \leq K \nu(B(0,1/2)^c)\, r^{d-\alpha}\,.
$$\qed

Theorem \ref{t:V} yields that $V$ is continuous on $\R\setminus\{ 0
\}$. Consequently, $V(x)\approx |x|^{\alpha-d}$ and $G(x,y)$ is
continuous on $B\times B\setminus\{ (x,y):x=y\}$.
\begin{lemat}\label{Greenwsrodku} $
    G(x,v) \approx |v-x|^{\alpha-d},\,\quad \mbox{ if }\quad |x|<1/2\,,\; |v|<3/4\,.
  $
\end{lemat}
We skip the proof as it is the same as the one of Lemma 6 in \cite{BS}.

We note that $\lim\limits_{x\to z} G(x,v)=0$ for every
$v\in B(0,1)$ and every point $z\in \sfera$ because 
the measures $\omega^x_{B(0,1)}$ weakly converge to $\delta_z$.
This is related to the regularity of $B(0,1)$, and it follows, e.g., 
from the estimate $$\omega^x_{B(x, 1-|x|)}(B(0,1)^c)\geq c\,,$$ which 
is a consequence of scaling, nondegeneracy of $\nu$ (compare
(\ref{Lduza})), and (\ref{IW}). 

We will employ the operator
 $$
  \Dynkinr{\phi}{x} = \frac{E^x\phi(X_{\tau_{B(x,r)}}) - \phi(x) } {E^x \tau_{B(x,r)}}\,,
 $$
  whenever the expression is well defined for given $\phi$, $r>0$
  and $x$. We note that ${\mathcal U}_r$ is implicitly used in
  \cite[Chapter III \S 17]{BgFt}.
Clearly, if $h$ is harmonic in $D$, $x\in D$, and
$r<\dist(x,D^c)$, then $\Dynkinr{h}{x}=0$. We note that
$$
  \Dynkin{\phi}{x}=\lim_{r\downarrow 0} \Dynkinr{\phi}{x}
$$
is the Dynkin characteristic operator, which was used
in \cite{BS} in a similar way.

We record the following observation (maximum principle).
\begin{lemat}\label{zasadaminimum}
  If there is $r>0$ such that $\Dynkinr{h}{x} > 0$ then $h(x)<\sup\limits_{y\in\R} h(y)$.
\end{lemat}

\begin{lemat}\label{Green<s}
  There exists $C=C(\alpha,\mu)$ such
  that
  $$
   G(x,v) < C s(v),\quad |x|<1/2\,,\quad 3/4 < |v| < 1\,.
  $$
\end{lemat}
\dowod 
By the strong Markov property we have
\begin{eqnarray*}
  s(v)
  &   =   & E^v \tau_B = E^v(\tau_A+\tau_{B(0,1)} \circ \theta_{\tau_A})
     =    E^v\tau_A + E^v E^{X_{\tau_A}} \tau_{B(0,1)} \\
  &   =   & E^v \tau_A + E^v s(X_{\tau_A})\,,\quad v\in \R\,,\quad A\subset {B(0,1)},
\end{eqnarray*}
which yields $\Dynkinr{s}{v}=-1$ for $v\in {B(0,1)}$ and $r<1-|v|$.

For $n\in \{1,2,\ldots\}$ and  $x\in B(0,1/2)$ we let 
$g(v) =G(x,v)$ and $g_n(v) =\min(G(x,v),n)$.
For $v \in B(x,1/8)^c$ we have that
$G(x,v) \leq c_1|x-v|^{\alpha-d}$ hence
$g_n(v)=G(x,v)$ provided $n \geq c_1 8^{d-\alpha}$. 
By harmonicity of $g$ on ${B(0,1)}\setminus\{x \}$, 
scaling property, (\ref{IW}) and (\ref{RK}) 
we obtain  that for $v\in {B(0,1)}\setminus B(0,3/4)$ and $r<\min(1-|v|,1/16)$ it holds
\begin{eqnarray*}
  \Dynkinr{g_n}{v} 
 &  =    &    \Dynkinr{(g_n-g)}{v} \\
 &  =   & \frac{1}{E^0\tau_{B(0,1)}}
           \int\limits_{B(0,1)} G(0,w)
           \int(g_n-g)(v+rw+z) \nu(dz)dw \\
  & \geq & \frac{-c_2}{s(0)} \int\limits_{B(0,1)} G(0,w)
           \int\limits_{B(x-v-rw,1/8)}|x-v-rw-z|^{\alpha-d} \nu(dz)dw \\
  & \geq & \frac{-c_2K}{s(0)} \int\limits_{B(0,1)} G(0,w)
           \nu(B(x-v-rw,1/8))dw \geq -c_3\,.
\end{eqnarray*}
If $a>c_3$ then
$$
\Dynkinr{(a s - g_n)}{v} = -a - \Dynkinr{g_n}{v} \leq -a + c_3 < 0.
$$
By scaling
\begin{eqnarray}\label{sdodatnie}
s(v) 
& \geq & E^v \tau_{B(v,1-|v|)}=(1-|v|)^\alpha E^0\tau_{B(0,1)} \\
& \geq & 4^{-\alpha} E^0\tau_{B(0,1)}\,, \quad |v|<3/4\,. \nonumber
\end{eqnarray}
Since $g_n(v) \leq n$, we see that $as(v) - g_n(v) > 0$ for $v\in B(0,3/4)$ provided $a >
n/(4^{-\alpha}E^0\tau_{B(0,1)} )$.

Let $a_0= \max[c_3,n/(4^{-\alpha}E^0\tau_{B(0,1)})] + 1$ and $h(v) =
a_0 s(v) - g_n(v)$.  We have $h(v)\geq 0 $ for $v\in
\overline{B(0,3/4)}$, $h(v)=0$ for $v\in {B(0,1)}^c$ and $\Dynkinr{h}{v} < 0$
for $v\in {B(0,1)}\setminus B(0,3/4)$, $r<\min(1-|v|,1/16)$.
Lemma \ref{zasadaminimum} and 
continuity of $h$ yields $h(v)\geq 0 $ in ${B(0,1)}$.  Since $g_n=g$ on
$B(0,3/4)^c$, the lemma follows.\qed
\vspace{2mm}

Lemma \ref{Green<s} and \ref{ET_est} yield the following conclusion:
\begin{equation}\label{oszacowanieG}
    G(x,v) \leq C (1-|v|)^{\alpha/2}\,,\quad |x|<1/2\,,\;3/4<|v|<1\,.
\end{equation}

\begin{lemat}\label{Green>s}
  There is
  $C=C(\alpha,\mu)$ such that
  $G(x,v) \geq C s(v)$ provided $|x|<1/2$ and $|v|<1$.
\end{lemat}

\dowod Let $x\in B(0,1/2)$. We fix $\varepsilon$ such that
(\ref{Lduza}) is satisfied.  Lemma \ref{Greenwsrodku} yields that
$G(x,v) \geq c_1>0$ for $v\in B(0,1-\varepsilon)$.  Let $n\in \{1,2,\ldots\}$
be such that $c_1\geq 2/n$.  By (\ref{oszacowanieG}) there is
$\eta>0$ such that $G(x,v) \leq 1/n$ for $v \in {B(0,1)}\setminus
B(0,1-\eta)$. Let $g(v) = G(x,v)$ and $g_n(v) =\min( g(v),1/n)$. We
have
$$
g_n(v) = g(v),\quad v\in {B(0,1)}\setminus B(0,1-\eta),
$$
and
$$
g(v) - g_n(v) \geq 2/n-1/n =1/n\,,\quad v\in B(0,1-\varepsilon)\,,
$$
hence by Lemma \ref{czescDynkina} for
$v\in {B(0,1)}\setminus \overline{B(0,1-\eta)}$ and $r<\min(1-|v|,(\varepsilon-\eta)/2)$ we obtain
\begin{eqnarray*}
  \Dynkinr{g_n}{v} 
  &   =   &  \Dynkinr{(g_n-g)}{v} \\
  &   =   & \frac{1}{s(0)}\int\limits_{B(0,1)} G_{B(0,1)}(0,w)
            \int(g_n-g)(v+rw+z) \nu(dz) dw\\
  & \leq & -\frac{1}{n} \frac{1}{s(0)}\int\limits_{B(0,1)} G_{B(0,1)}(0,w)
           \nu(B(v+rw,1-\varepsilon)) dw
            \leq -\frac{c_2}{n}.
\end{eqnarray*}
For $a>0$ we have
$$
\Dynkinr{(a g_n - s)}{v} \leq -c_2a/n + 1\,,\quad v\in
{B(0,1)}\setminus B(0,1-\eta)\,.
$$
This is negative if $a > n/c_2$.  Furthermore $s(v)\leq
c_3$ for $v\in {B(0,1)}$ and $g_n(v) \geq c_4>0$ for $v\in B(0,1-\eta)$. Thus
$a g_n(v) - s(v) \geq ac_4 - c_3 > 0$ for $v\in B(0,1-\eta)$ if only
$a > c_3/c_4$.  Note that our estimates do not depend on $x$, provided
$|x|<1/2$.  Let $a_0 = \max(c_3/c_4,n/c_2) +1$ and $h(v)
= a_0 g_n(v) - s(v)$.  We have $h(v) \geq 0$ for $v\in \overline{
  B(0,1-\eta)}$ and $\Dynkinr{h}{v} <0 $ for $v\in {B(0,1)}\setminus
B(0,1-\eta)$. By Lemma \ref{zasadaminimum} and the continuity of $h$
we get $ h(v) \geq 0 $ in ${B(0,1)}$ and the lemma follows.
\qed

\vspace{15pt}
\noindent
{\em Proof of Proposition \ref{OG}}.\/ The estimate is a consequence of (\ref{sdodatnie}),
Lemma~\ref{Greenwsrodku}, \ref{Green<s}, and \ref{Green>s}.  \qed

Maciej Lewandowski \cite{Lewandowski} has informed us that he recently proved the converse of
the inequality in Lemma \ref{ET_est}. This implies
\begin{equation}\label{e:gL}
    G(x,v) \approx (1-|v|^2)^{\alpha/2} |v-x|^{\alpha-d}\,,\quad |x|<1/2\,,\; |v|<1\,.
\end{equation}
We will not use (\ref{e:gL}) in the sequel; the less explicit estimate in Lemma \ref{Green>s} suffices for our purposes.
Note that the asymptotic of $G$ at the pole is different when
$d=1\leq \alpha$, see, e.g., \cite{BB2}.

\begin{lemat}\label{lem_Harnack}
(RK) implies Harnack's  inequality for all $d\in \{1,2,\ldots\}$ and
 $\alpha\in (0,2)$.
\end{lemat}
\dowod 
By translation and scaling invariance of the  class of
harmonic functions and by a
covering argument we only need to verify that
\begin{displaymath}
  u(0)\leq c\, u(x)\,,\quad |x|<1/2\,,
\end{displaymath}
whenever $u$ is nonnegative on $\R$ and regular harmonic on $B(0,1)$.
For this to hold it is sufficient to have, with the same constant $c$,
\begin{equation}\label{e:pp}
  P(0,y)\leq c\, P(x,y)\,,\quad |x|<1/2\,,\quad |y|>1\,.
\end{equation}
If $d=1$, (\ref{e:pp}) follows from (\ref{e:jp}).
Thus we only need to examine the case $d>\alpha$.
By the decomposition ${B(0,1)}=B(0,1/2)\cup [{B(0,1)}\setminus
B(0,1/2)]$, (\ref{e:P}), Proposition~\ref{OG}, (\ref{RK}), and the fact that $s$ is bounded away
from zero on compact subsets of ${B(0,1)}$ (comp. (\ref{sdodatnie})), we obtain
\begin{eqnarray*}
  P(0,y)&\approx&\int_{B(y,1)} s(y-v) |y-v|^{\alpha-d}\nu(dv)
  \approx\int_{B(y,1)} s(y-v) \nu(dv)\\
  &\leq& c\,\int_{B(y,1)} s(y-v)|y-v-x|^{\alpha-d} \nu(dv) \\
  & \approx & P(x,y)\,, \quad |x|<1/2\,,\; |y|>1\,.
\end{eqnarray*}
\qed

\vspace{15pt}
{\it Proof of Theorem \ref{Harnack}}\/.
See Lemma \ref{n_RK} and Lemma \ref{lem_Harnack}.
\qed

\vspace{5mm}
We conclude with a few remarks and open problems.

By translation and dilation invariance of the class of considered harmonic
functions, and by a covering argument Harnack's  inequality
holds for every compact subset of every connected domain of harmonicity.
We note that : (1) it does not generally hold for disconnected open sets, as the support 
of $y\mapsto P(x,y)$ may be smaller than $B(0,1)^c$ (see \ref{IW}), (2) it does hold for all open sets if $\nu$ is isotropic 
(this follows from (\ref{e:jp}), or see \cite{BBsm1999}). 

We consider the following examples of measures $\nu$. 
(RK) holds for $\nu_1(dy)\approx |y|^{-d-\alpha}dy$ (both sides of (\ref{RK}) may be explicitly estimated).
Next, let $\xi\in\sfera$, $0<r<\sqrt{2}$, and $C=\sfera\cap[B(\xi,r)\cup B(-\xi,r)]$. 
(RK) holds for $\nu_2(dy)=\indyk{C}(y/|y|)|y|^{-d-\alpha}dy$, see \cite{BS}. 

On the other hand, consider balls $B_{n}\subset B'_{n}$ centered at $\sfera$, 
with radii $4^{-n}$ and $2^{-n}$, respectively, and such that $\{B'_n\}$ are pairwise disjoint.
Let $C=\bigcup_{n\geq n_0}B_n$ and let $\nu_3(dy)=\indyk{C}(y/|y|)|y|^{-d-\alpha}dy$.
If $d-1>\alpha$ then (RK) does not hold for $\nu_3$ (\cite{BS}) even though it 
is bounded by $\nu_1$.

Let $B_{\xi,r}=B(\xi,r)\cap \sfera$. By integrating in polar coordinates we can 
give the characterization of relative Kato condition in terms of its spectral measure $\mu$
and $B_{\xi,r}$ (comp. \cite{BS}). 
Let $d-\alpha>1$. (RK) holds for $\nu$ if
and only if
\begin{equation}
  \label{eq:rkr2}
  \int_{B_{\xi,r}}(|\eta-\xi|/r)^{\alpha-(d-1)}\mu(d\eta)
\leq c\, 
  \mu(B_{\xi,r})\,,\quad \xi \in \sfera\,,\;0<r<c\,.
\end{equation}
In the case $d=2$, $\alpha=1$, (RK) is equivalent to
\begin{equation}
  \label{eq:rkr3}
  \int_{B_{\xi,r}}\log(2r/|\eta-\xi|)\mu(d\eta)
\leq c\, 
  \mu(B_{\xi,r})\,,\quad \xi \in \sfera\,,\;0<r<c\,.
\end{equation}
In the case of $d=2$ and $\alpha>1$ (RK) is always satisfied. We omit the proofs.
\begin{wniosek}
If $d-1<\alpha$ then
Harnack's inequality holds for $\gener$.
\end{wniosek}
This may be extended as follows.
We will say that $\nu$ is a {\it strict} $\gamma$-measure if
\begin{equation}\label{s_nu_gamma} 
  \nu(B(x,r)) \approx r^{\gamma}\,,\quad \mbox{ provided } x\in \supp \nu\,,\; |x|=1\,,\; 0<r<1/2\,,
\end{equation}
compare (\ref{nu_gamma}).
Of course, if $\nu$ is a (strict) $\gamma$ measure on $\sfera$ than $\mu$ is a 
(strict) $(\gamma-1)$-measure (on $\sfera$).  
This observation and (\ref{eq:rkr2}) yield the following conclusion, which we state without proof.
\begin{wniosek}
If $\nu$ is a strict $\gamma$-measure with $\gamma>d-\alpha$, then Harnack's inequality holds for $\gener$.
\end{wniosek}
The example of $\nu_3$ shows the importance of the strictness assumption.
We interpret (RK) as a property of balance or firmness of $\nu$. 
As such it is close to the reverse H\"older condition with exponent $q>d/\alpha$, see \cite{BS}.

If $\mu({\xi})>0$ for some $\xi\in\sfera$ then $\nu$ is a $1$-measure only. 
By Theorem~\ref{t:nu} the potential kernel $V$ is unbounded on $\sfera$ if $1>d-2\alpha$ 
(in fact, if $1\geq d-2\alpha$, see \cite[Theorem 1.1]{W}, \cite{BS}).
That $V$ may be infinite on rays emanating from the origin 
shows that harmonic functions cannot be defined pointwise by means of
$\gener$. In general they even lack finiteness in the domain of
harmonicity (but see \cite{BBsm1999} and \cite{PiSa} in this connection).
Thus the potential-theoretic properties of the operators $\gener$ are very diverse among considered measures $\nu$.
This is in sharp contrast with the fact that the exponents $\Phi$ (see (\ref{Phi})) 
are all comparable and the same is true of the corresponding Dirichlet forms (\cite{FOT}, see also
\cite{D}). The boundary potential theory of $\gener$ will generally be 
very different from that of the fractional Laplacian 
(see \cite[p. 199]{SztonykCOLL} for a simple remark on this subject).

We like to mention a number of further interesting problems and references: (1) characterization of
continuity and higher order regularity of $V$ on $\sfera$ (\cite{BsLn}), 
(2) the boundary Harnack principle (comp. \cite{Bo1, BSS1, SztonykCOLL}), the corresponding approximate 
factorization of $G(x,v)$ for {\it all}\/ $x,v\in B(0,1)$  (comp. \cite{Bo2, Ha, CK} and Proposition \ref{OG} above), 
and related boundary problems (comp. \cite{BsBn, BrKl}),
(3) study of other L\'evy measures which are in the form of a product in
polar coordinates, (4) study of similar nonlocal operators $\gener$ 
which are not translation invariant (\cite{BrBsKs, SV}).

{\small K. Bogdan (corresponding author), 
Institute of Mathematics, Polish Academy of Sciences,
Institute of Mathematics, Wroc\l{}aw University of Technology,

{\it E-mail address: bogdan@pwr.wroc.pl}

\vspace{2mm}
P. Sztonyk, Institute of Mathematics, Wroc\l{}aw University of Technology,
Wybrze{\.z}e Wyspia\'nskiego 27, 50--370 Wroc\l{}aw, Poland

{\it E-mail address: sztonyk@im.pwr.wroc.pl}}


\begin{thebibliography}{2}

\bibitem{BsBn} R. Ba\~nuelos, K. Bogdan {\it Symmetric stable processes in cones},
 Potential Anal. 21 (2004), no. 3, 263--288.

\bibitem{Bs1}R. F. Bass {\it Diffusions and elliptic operators}, Probability
and its Applications (New York). Springer-Verlag, New York, 1998.

\bibitem{Bs2} R. F. Bass, {\it Stochastic differential equations
  with jumps},  Probab. Surv.  1  (2004), 1--19 (electronic).

\bibitem{BsKn}
R. F. Bass, M. Kassmann, {\it Harnack inequalities for non-local
  operators of variable order}, Trans. Amer. Math. Soc. 357 (2005), 837-850.

\bibitem{BrBsKs}
M. Barlow, R. F. Bass, M. Kassmann, 2005, preprint.

\bibitem{BsLn} R. F. Bass, D. A. Levin, {\it Harnack inequalities
    for jump processes,} Potential Anal. 17(4)(2002), 375-388.

\bibitem{BgFt} C.~Berg and G.~Forst, {\it Potential Theory on 
Locally Compact Abelian Groups}. Springer-Verlag, 1975.

\bibitem{BrKl} K. Burdzy,T. Kulczycki {\it Stable processes have thorns}, 
Ann. Probab. 31 (2003), no. 1, 170--194.

\bibitem{Br} J. Bertoin, {\it L\' evy processes,} Cambridge
  University Press, Cambridge, 1996.

\bibitem{BH} J. Bliedtner, W. Hansen, 
{\it Potential theory. An analytic and probabilistic approach to balayage} 
Springer-Verlag, Berlin, 1986.

\bibitem{BGB} R. M. Blumenthal and R. K. Getoor, {\it Markov Processes
    and Potential Theory,} Pure Appl. Math., Academic Press Inc., New
  York 1968.

\bibitem{BGR} R. M. Blumenthal, R.K. Getoor, and D.B. Ray, {\it On the
    distribution of first hits for the symmetric stable processes,}
  Trans. Amer. Math. Soc. 99 (1961), 540-554.

\bibitem{Bo1} K. Bogdan, {\it The boundary Harnack principle for the fractional Laplacian}, 
Studia Math. 123 (1997), no. 1, 43--80.
 
\bibitem{Bo2} K. Bogdan, {\it Sharp estimates for the Green
    function in Lipschitz domains,}
  J. Math. Anal. Appl. 243 (2000), no. 2, 326--337.

\bibitem{BBsm1999}
K. Bogdan, T. Byczkowski, {\it Potential theory for the
  {$\alpha$}-stable {S}chr\"odinger operator on bounded {L}ipschitz domains,}
Studia Math. 133 (1999) no.~1, 53--92.

\bibitem{BB2} K. Bogdan, T. Byczkowski, {\it Potential theory of the
    Schr\"odinger operator based on the fractional Laplacian}, Prob.
  Math. Stat. 20 (2000), 293-335.

\bibitem{BSS1} K. Bogdan, A. St\'os, P. Sztonyk, {\it Potential theory
    for L\'evy stable processes,} Bull. Polish. Acad. Sci. Math. 50(3)
  2002, 361--372.

\bibitem{BS} K. Bogdan, P. Sztonyk, {\it Harnack's  inequality for
    stable L{\'e}vy processes,} Potential Analysis 22 (2) (2005), 1333--150.

\bibitem{BZ} K. Bogdan, T. {\.Z}ak, {\it On {K}elvin
    transformation} (2004), to appear in J. Theoretical Probability

\bibitem{CS}
Z.-Q.~Chen and R.~Song, Martin boundary and integral representation for
harmonic functions of symmetric stable processes. {\it J. Funct. Anal.} 
{\bf 159} (1998),  267--294.

\bibitem{CK} Z.-Q. Chen, P. Kim, 
{\it Green function estimate for censored stable processes},
Probab. Theory Related Fields 124 (2002), no. 4, 595--610.

\bibitem{Ch} K. L. Chung, 
{\it Lectures from Markov processes to Brownian motion}, 
Springer-Verlag, New York-Berlin, 1982. 

\bibitem{Ch1} K. L. Chung, {\it Doubly-Feller process with multiplicative functional},
Seminar on stochastic processes, 1985 (Gainesville, Fla., 1985), 63--78, Progr. Probab. Statist., 12, 
Birkhäuser Boston, Boston, MA, 1986. 

\bibitem{ChZ} K. L. Chung, Z. Zhao, {\it From Brownian motion to
   Schr\"odinger's equation,} Springer - Verlag, New York, 1995.

\bibitem{D} B. Dyda, {\it On comparability of integral forms} (2004), to appear in J. Math. Anal. Appl.

\bibitem{Dz} J. Dziuba\'nski, Asymptotic behaviour of
  densities of stable semigroups of measures, Probab. Theory Related
  Fields 87 (1991), 459-467.

\bibitem{FOT}
Masatoshi Fukushima, Y{\=o}ichi {\=O}shima, and Masayoshi Takeda.
{\it Dirichlet forms and symmetric {M}arkov processes}, Walter de Gruyter \& Co., Berlin, 1994.

\bibitem{G} R. K. Getoor, {\it First passage times for symmetric
    stable processes in space,} Trans. Amer. Math. Soc. 101 (1961),
  75-90.

\bibitem{GT}
D. Gilbarg, N. S. Trudinger, 
{\it Elliptic partial differential equations of second order},
Springer, Berlin, 2001.

\bibitem{GH} P. G{\l}owacki, W. Hebisch, {\it Pointwise estimates for
    densities of stable semigroups of measures,} Studia Math. 104
  (1993), 243-258.

\bibitem{Ha} W. Hansen, {\it Uniform boundary Harnack principle and generalized triangle property} (2004), 
preprint


\bibitem{Hi} S. Hiraba {\it Asymptotic estimates for densities of
 multi-dimensional stable distributions},  Tsukuba J. Math.  27
 (2003),  no. 2, 261--287.

\bibitem{Ho}
W. Hoh {\it Pseudo differential operators generating Markov
  processes}, Habilitationsschrift, Universit\"at Bielefeld 1998.

\bibitem{IW} N. Ikeda, S. Watanabe, {\it On some relations between the
    harmonic measure and the L\'evy measure for a certain class of
    Markov processes,} J. Math. Kyoto Univ. 2-1 (1962), 79-95.

\bibitem{Jc1}
N. Jacob, {\it Pseudo differential operators and Markov
  processes. Vol. I. Fourier analysis and semigroups},
Imp. Coll. Press, London, 2001.

\bibitem{Jc2}
 N. Jacob,
 {\it Pseudo-Differential Operators and Markov Processes,
 Vol. II : Generators and Their Potential Theory,}
 Imperial College Press, London, 2002.

\bibitem{JS} N. Jacob, R. Schilling {\it L\'evy-type processes 
and pseudodifferential operators}.  L\'evy processes,  139--168, 
Birkhauser Boston, Boston, MA, 2001.

\bibitem{Jk} T. Jakubowski, {\it The estimates for the Green function
    in Lipschitz domains for the symmetric stable processes,} (2002),
  to appear in Colloq. Math.

\bibitem{K1}
T{.} Kulczycki,
\emph{Properties of Green function of symmetric stable process},
Probab{.} Math{.} Statist{.} 17(2) (1997), 339--364.

\bibitem{La}
N. S. Landkof {\it Foundations of Modern Potential Theory}.
Springer-Verlag, New York, 1972.

\bibitem{Lewandowski} M. Lewandowski, private communication (2005).

\bibitem{Millar} P. W. Millar, {\it First passage distributions of
    processes with independent increments,} Ann. Probab., 3 (1975) no.
  2, 215-233.

\bibitem{Pi} J. Picard, 
{\it Density in small time at accessible points for jump processes,} 
Stochastic Process. Appl. 67 (1997), no.~2, 251--279.

\bibitem{PiSa}
J. Picard, C. Savona, 
{\it Smoothness of harmonic functions for processes with jumps,}
Stochastic Process. Appl. 87 (2000), no.~1, 69--91.

\bibitem{RSV}
M. Rao, R. Song and Z. Vondra\v{c}ek, 
{\it Green function estimates and Harnack inequality 
for subordinate Brownian motions} (2004), to appear in Potential Analysis.

\bibitem{Sato} K.-I. Sato, {\it L\' evy Processes and Infinitely Divisible 
             Distributions}, Cambridge University Press, 1999.

\bibitem{SV} 
R. Song, Z. Vondra\v{c}ek, {\it Harnack inequality
    for some classes of Markov processes}, Math. Z.  246  (2004),
  no. 1-2, 177--202.

\bibitem{S1} D. W. Stroock {\it Diffusion processes associated with L\'evy generators},
Z. Wahrscheinlichkeitstheorie und Verw. Gebiete 32 (1975), no. 3,
209--244.

\bibitem{S2} D. W. Stroock, {\it Markov processes from K. It\^o's perspective}, Ann. of Math. Stud., 155, Princeton Univ. Press, 
Princeton, NJ, 2003

\bibitem{Sztonyk} P. Sztonyk, {\it On harmonic measures for L\'evy
    processes,} Prob. Math. Statist. 20(2) (2000), 383-390.

\bibitem{SztonykCOLL} P. Sztonyk, {\it Boundary potential theory for
    stable processes,} Colloq. Math., 95 (2003) no. 2, 191-206.

\bibitem{Taylor} S. J. Taylor, {\it Sample path properties of a
    transient stable process,} J. Math. Mech. 16 (1967), 1229-1246.

\bibitem{W} T. Watanabe {\it Asymptotic estimates of
multi-dimensional stable densities and their applications},
preprint (2004).

\bibitem{Wu}
J.-M. Wu {\it Harmonic measures for symmetric stable processes},
Studia Math.  149  (2002),  no. 3, 281--293.

\end{thebibliography}
\end{document}